\documentclass[10pt]{article}
\begin{document}

\title{Can We Prove Goldbach's Conjecture?}
\author{Danilo Mauro\\
Via Bancaria 2\\
34071 Cormons - Italy\\
danilo.mauro@istruzione.it}
\date{}
\maketitle

\noindent In this paper we will propose a strategy to prove Goldbach's conjecture: 
every even integer greater than 2 can be written as the sum of two primes.

\section{A set of conditions for Goldbach's partitions}

Let us start this paper by considering a generic even number $q$ and its partitions as the sum of two integers $n_1$ and $n_2$ greater than 1, i.e.:
\begin{equation}
q=n_1+n_2. \label{zero}
\end{equation}
In order to have a Goldbach's partition \cite{guy} of $q$, the numbers $n_1$ and $n_2$ must be prime numbers, so they have to satisfy both the following conditions:
\begin{itemize}
\item the remainders of the ratio of $n_1$ and all the primes smaller than $\sqrt{q}$:
\begin{equation}
p_0=2,  \;p_1=3, \;p_2=5, \;\dots, \; p_j \label{uno}
\end{equation}
must be different than zero, unless $n_1$ itself coincides with a prime in the list (\ref{uno}). In fact, if the remainder of the ratio of $n_1$ and one or more primes $p_i$ in (\ref{uno}) were zero then the integer $n_1$ would not be a prime number. We can limit ourselves to consider only the primes smaller than $\sqrt{q}$ because if $n_1$ is composite then it must appear at least one factor smaller than $\sqrt{q}$ in its factorization as a product of prime numbers.
\item unless $n_2$ coincides with a prime in the list (\ref{uno}), the remainder of the ratio of $n_1$ and a generic prime $p_i$ in (\ref{uno}) must be different than the remainder of the ratio of the even number $q$ and the same prime $p_i$. In fact, this implies that the remainder of the ratio of the second number $n_2$ in the partition of $q$ and all the primes in (\ref{uno}) is different than zero. Therefore $n_2$ cannot be a multiple of $2, 3,  \dots, p_j$, consequently it is a prime.\end{itemize}

The previous rules can be applied to the very simple case in which we consider only the prime number 2. The first rule states that the number $n_1$ cannot be a multiple of 2, so it must be an odd number. In this case the second rule is equivalent to the first one, and it implies that also $n_2$ must be an odd number. Of course there is only one exception given by $4=2+2$ in which both $n_1$ and $n_2$ coincide with 2 itself. 

\section{Bounds in the number of Goldbach's partitions}

In the previous section we have found out that, except for the case $4=2+2$, both the numbers $n_1$ and $n_2$, which appear in a  generic partition (\ref{zero}) of $q$, must be odd numbers in order to have a Goldbach's partition. 
The total number of ways in which we can rewrite  $q$ as the sum of two odd numbers greater than 1 (which is not a prime number) is $n=q/2-2$. For example if $q=6$ there is only $n=1$ partition $6=3+3$, if $q=8$ there are $n=2$ symmetric partitions $8=5+3=3+5$, and so on. 

According to the rules given in the previous section, for each prime number $p_i$ with $1\le i\le j$ we have to disregard two bad remainders: 0 (to be sure that $n_1$ is a prime) and the remainder of the ratio of $q$ and $p_i$ (to be sure  that $n_2$ is a prime)\footnote{In this way we disregard Goldbach's partitions in which $n_1$ or $n_2$ is a prime $p_i$ in the list (\ref{uno}) but our only goal is to prove that every even number admits {\it at least one} Goldbach's partition.}. Consequently in the worst case, i.e. when the two bad remainders of $p_i$ do not coincide, there will be only one good remainder of 3, three good remainders of 5, in general $p_i-2$ good remainders of the prime number $p_i$. 

For every even number $q \ge 6$, we must take into account the bad remainders of the prime number $p_1=3$. When we change the number $n_1$, all the remainders of the ratio of $n_1$ and 3 are equally distributed, see Table 1.
\begin{table}[htdp]
\caption{Remainders of the ratio of $n_1$ and the first prime numbers}
\begin{center}
\begin{tabular}{|c|ccccccccc|}
$n_1$ & r(3) & r(5) & r(7) & r(11) & r(13) & r(17) & r(19) & r(23) & r(29) \\ \hline
3 & 0 & 3 & 3 & 3 & 3 & 3 & 3 & 3 & 3 \\
5 & 2 & 0 & 5 & 5 & 5 & 5 & 5 & 5 & 5 \\
7 & 1 & 2 & 0 & 7 & 7 & 7 & 7 & 7 & 7 \\
9 & 0 & 4 & 2 & 9 & 9 & 9 & 9 & 9 & 9 \\
11 & 2 & 1 & 4 & 0 & 11 & 11 & 11 & 11 & 11 \\
13 & 1 & 3 & 6 & 2 & 0 & 13 & 13 & 13 & 13 \\
15 & 0 & 0 & 1 & 4 & 2 & 15 & 15 & 15 & 15 \\
17 & 2 & 2 & 3 & 6 & 4 & 0 & 17 & 17 & 17 \\
19 & 1 & 4 & 5 & 8 & 6 & 2 & 0 & 19 & 19 \\
21 & 0 & 1 & 0 & 10 & 8 & 4 & 2 & 21 & 21 \\
23 & 2 & 3 & 2 & 1 & 10 & 6 & 4 & 0 & 23 \\
25 & 1 & 0 & 4 & 3 & 12 & 8 & 6 & 2 & 25 \\
27 & 0 & 2 & 6 & 5 & 1 & 10 & 8 & 4 & 27 \\
29 & 2 & 4 & 1 & 7 & 3 & 12 & 10 & 6 & 0 \\
31 & 1 & 1 & 3 & 9 & 5 & 14 & 12 & 8 & 2 \\
\end{tabular}
\end{center}
\label{default}
\end{table}%
The worst situation happens when the first 2 remainders of the ratio of $n_1$ and 3 are just the bad remainders\footnote{For example, if $q=20$ then the bad remainders of the ratio of $n_1$ and 3 are 0 (in a Goldbach's partition $n_1$ cannot be a multiple of 3) and 2 (in a Goldbach's partition $n_2=q-n_1$ cannot be a multiple of 3) but 0 and 2 are just the two numbers which appear at the beginning of the column r(3) in Table 1.} of 3. After having disregarded the bad remainders of 3 we have at least 
\begin{equation}
A=\frac{1}{3}\cdot (n-2) = \frac{1}{3}\cdot n-\frac{2}{3}>\frac{1}{3}\cdot n -1 \label{agg1}
\end{equation}
possible partitions for the number $q$. 

When we consider even numbers $q\ge 28$, 
we have to be sure that there is no problem with both 3 and 5. Also the remainders of $n_1$ and 5 are equally distributed, see Table 1. Furthermore, when we change the number $n_1$, which appears in the partition of $q$, we have that the remainders of 3 and 5 are independent and every remainder of 3 is associated with every remainder of 5. Also in this case the worst situation happens when the first two remainders of 5, associated with the only good remainder of 3, are bad remainders of 5. This implies that the number of good partitions, after having excluded the bad remainders of both 3 and 5, is at least equal to
\begin{equation}
A = \frac{3}{5}\cdot \biggl(\frac{1}{3}\cdot n -1\biggr)-2>\frac{1}{3}\cdot\frac{3}{5}\cdot n
-\frac{3}{5}-2>   \frac{1}{3}\cdot\frac{3}{5} \cdot n-3.   \label{agg2}
\end{equation}

\section{Towards a proof of Goldbach's conjecture}

If the number $q$ satisfies $p_j^2+3\le q\le p_{j+1}^2+3$ then we have to take into account all the primes up to $p_j$. All the remainders of $n_1$ and all the primes $p_1, \dots, p_j$ are equally distributed and independent. Let us make the following conjecture: the total number of Goldbach's partition for $q$ is greater than
\begin{equation}
A= \frac{1}{3}\cdot\frac{3}{5}\cdot\frac{5}{7}\cdot\frac{9}{11}\cdots \frac{p_j-2}{p_j}\cdot n-(p_j-2), \label{five}
\end{equation}
where $p_j$ is the greatest prime number smaller than $\sqrt{q}$. The term $(p_j-2)$ in (\ref{five}) plays the role of the term $-1$ in (\ref{agg1}) and of the term $-3$ in (\ref{agg2}) and corresponds to the worst situation among all the possible ones.

The local minima $A_m$ of this function correspond to the points in which a new prime starts to be taken into account, i.e. when $q=p_m^2+3$ for a certain prime $p_m$. The number of partitions $n=q/2-2$ becomes:
\begin{displaymath}
n=\frac{p^2_m+3}{2}-2=\frac{p^2_m-1}{2}.
\end{displaymath}
Replacing the previous expression in (\ref{five}) we get:
\begin{displaymath}
A_m= \frac{1}{3}\cdot\frac{3}{5}\cdot\frac{5}{7}\cdot\frac{9}{11}\cdots \frac{p_m-2}{p_m}\cdot \frac{p_m^2-1}{2}-(p_m-2).
\end{displaymath}
It is easy to realize that, for $m\ge 11$, $A_m$ can be rewritten as follows:
\begin{equation}
A_m=\frac{9}{7}\cdot \frac{15}{13}\cdot \frac{21}{19}\cdot \frac{27}{23}\cdot\frac{35}{31}\cdots \frac{p_m^2-1}{2\cdot p_m}-(p_m-2). \label{3bis}
\end{equation}
The product of the first five terms, which are always present for $m\ge 11$, is:
\begin{displaymath}
\frac{9}{7}\cdot \frac{15}{13}\cdot \frac{21}{19}\cdot \frac{27}{23}\cdot\frac{35}{31}>2.
\end{displaymath}
When $m>11$ there will be also other numerical fractions on the RHS of Eq. (\ref{3bis}), e.g. $39/37$, but all of them will be  greater than one. So we can conclude that, for every $m\ge 11$
\begin{displaymath}
A_m > \frac{p_m^2-1}{p_m}-(p_m-2) = 2-\frac{1}{p_m} > 1. 
\end{displaymath}
So, Goldbach's Conjecture is proven, provided that we succeed in giving a rigorous proof of Eq. (\ref{five}).

\section*{Acknowledgments}

I want to thank Robert Wilms for a counterexample that he found in the first version of this paper.

\end{document}